\documentclass[12pt]{article}

\usepackage{amsfonts}
\usepackage{amsmath}
\usepackage{amssymb}
\usepackage{color}
\usepackage{amsthm}
\usepackage{mathtools}
\usepackage{hyperref}

\newtheorem{theorem}{Theorem}
\newtheorem{corollary}[theorem]{Corollary}

\title{On Stein's method for stochastically monotone single-birth chains} \author{Fraser Daly\footnote{Department of Actuarial Mathematics and Statistics and the Maxwell Institute for Mathematical Sciences, Heriot-Watt University, Edinburgh EH14 4AS, UK.  E-mail: f.daly@hw.ac.uk}} \date{\today}

\begin{document}

\maketitle

\noindent{\bf Abstract} 
We discuss Stein's method for approximation by the stationary distribution of a single-birth Markov chain, in conjunction with stochastic monotonicity and similar assumptions. We use bounds on the increments of the solution of Poisson's equation for such a process. Applications include rates of convergence to stationarity, and bounding the total variation distance between the stationary distributions of two Markov chains in the case where one transition matrix dominates the other.

\vspace{12pt}


\vspace{12pt}


\section{Introduction}\label{sec:intro}

Let $\{Z_t:t=0,1,\ldots\}$ be a positive recurrent single-birth Markov chain in discrete time on the non-negative integers, with transition matrix $P$ whose $(i,j)$th entry we denote by $P_{i,j}$. Throughout we assume that $P_{i,i+1}>0$ for each $i\in\mathbb{Z}^+=\{0,1,2,\ldots\}$, and that $P_{i,j}=0$ for $j>i+1$. See Corollary 3.4 of \cite{l10} for conditions under which such a chain is positive recurrent. We let this chain have stationary distribution denoted by $(\pi_0,\pi_1,\pi_2\ldots)$, and $\pi$ be a random variable with this stationary distribution.

For a given function $h:\mathbb{Z}^+\to\mathbb{R}$, we let $f:\mathbb{Z}^+\to\mathbb{R}$ denote the solution to Poisson's equation
\begin{equation}\label{eq:poisson}
h(i)-\mathbb{E}h(\pi)=f(i)-\sum_{j=0}^\infty P_{i,j}f(j),
\end{equation}
with $f(0)=0$. In this note we will exploit this equation in conjunction with Stein's method to find explicit bounds in approximation by the stationary distribution of our single-birth chain. We will introduce the elements of Stein's method that we need in our work, but refer the reader to \cite{r11} and references therein for an introduction to this technique.

In this section we will derive Theorem \ref{thm:main} below, which gives an explicit bound in approximation by our stationary distribution, before exploring applications of this bound in Sections \ref{sec:convergence} and \ref{sec:comparison}. In  these applications it will be convenient to exploit stochastic monotonicity of the transition matrix $P$ or other similar monotonicity properties. Recall that $P$ is defined to be stochastically monotone if the distributions in successive rows of $P$ are stochastically non-decreasing; in this case stochastic ordering is preserved under transitions taken according to $P$. See \cite{d68} for further details.

One central aim of our work is to demonstrate how assumptions of stochastic monotonicity can be used in deriving bounds in approximations by the distribution of $\pi$ using Stein's method, in a similar spirit to \cite{dlu12}. This is somewhat different to earlier work on Stein's method for stationary distributions of Markov processes, including that of Brown and Xia \cite{bx01}, who considered approximation by the stationary distribution of a continuous-time birth-death process without any assumptions of monotonicity, and whose applications are of a quite different flavour to ours. In the discrete-time setting, \cite{bn19,rr19} have recently used Poisson's equation as a starting point for Stein's method, but in the case of a finite state space and again without monotonicity assumptions. The focus of our work is thus somewhat different to these other papers. While there are many other tools available in the literature for tackling examples and applications such as those we consider here (see the discussion and references in the examples below), our main purpose here is to show how Stein's method can be added to such a toolkit.  

Letting $\widehat{h}(j)=h(j)-\mathbb{E}h(\pi)$, Jiang \emph{et al$.$} \cite{jly14} have shown that Poisson's equation \eqref{eq:poisson} is solved by a function $f$ satisfying 
$f(j+1)-f(j)=-m_j(h)$ for $j\geq0$, where $m_0(h)=\widehat{h}(0)/P_{0,1}$ and
\begin{equation}\label{eq:jly}
m_j(h)=\frac{1}{P_{j,j+1}}\left(\widehat{h}(j)+\sum_{k=0}^{j-1}m_k(h)\sum_{l=0}^kP_{j,l}\right),
\end{equation}
for $j\geq1$; see their Theorem 2.1. Our approximation results will make use of bounds on $\sup_{j\in\mathbb{Z}^+}|m_j(h)|$ for $h\in\mathcal{H}=\{h:\mathbb{Z}^+\to\mathbb{R}:|h(j)|\leq1\mbox{ for all }j\}$. This is the relevant class of functions for us since we derive bounds on the total variation distance between $\pi$ and a non-negative, integer-valued random variable $X$. This total variation distance is defined by
\[
d_{TV}(X,\pi)=\sup_{h\in\mathcal{H}}|\mathbb{E}h(X)-\mathbb{E}h(\pi)|.
\] 
Following Stein's method, this will be bounded by rewriting the right-hand side using Poisson's equation \eqref{eq:poisson}: 
\begin{align*}
\mathbb{E}h(X)-\mathbb{E}h(\pi)&
=\mathbb{E}f(X)-\mathbb{E}\sum_{j=0}^\infty P_{X,j}f(j)\\
&=\sum_{j=0}^\infty f(j)\left[\mathbb{P}(X=j)-\sum_{i=0}^\infty\mathbb{P}(X=i)\mathbb{P}(Z_1=j|Z_0=i)\right]\\
&=\sum_{j=0}^\infty\Delta f(j)\left[\mathbb{P}(X>j)-\sum_{k=j+1}^\infty\sum_{i=0}^\infty\mathbb{P}(X=i)\mathbb{P}(Z_1=k|Z_0=i)\right],
\end{align*}
where we recall that $f(0)=0$ and write $\Delta f(j)=f(j+1)-f(j)$. This immediately yields the following result.
\begin{theorem}\label{thm:main}
Let $\{Z_t:t=0,1,\ldots\}$ be a positive recurrent single-birth Markov chain on $\mathbb{Z}^+$ with transition matrix $P$ and stationary distribution $(\pi_0,\pi_1,\ldots)$. Let $X$ be a random variable supported on $\mathbb{Z}^+$. Then
\[
d_{TV}(X,\pi)\leq\sup_{h\in\mathcal{H}}\sup_{l\in\mathbb{Z}^+}|m_l(h)|\sum_{j=0}^\infty\left|\mathbb{P}(X>j)-\sum_{k=j+1}^\infty\sum_{i=0}^\infty\mathbb{P}(X=i)P_{i,k}\right|.
\]
\end{theorem} 
Note that this result makes use of bounds on the increments of the solution $f$ of Poisson's equation, not bounds on $f$ itself. In the remainder of this section we note some cases in which a bound on the increments of $f$ may be easily found and which we will use as running examples to illustrate some applications of Theorem \ref{thm:main} in conjunction with assumptions of stochastic monotonicity and domination throughout the remainder of this note. In Section \ref{sec:convergence} we use Theorem \ref{thm:main} to establish rates of convergence to stationarity for stochastically monotone single-birth processes. Then, in Section \ref{sec:comparison} we consider the approximation by $\pi$ of the stationary distribution of a Markov chain whose transition matrix either dominates, or is dominated by, $P$. These applications will be illustrated using the examples we introduce here.  

\subsection{Example: Birth--death chain}\label{sec:bd}

Suppose that $\{Z_t:t=0,1,\ldots\}$ is a birth-death chain with $P_{i,i+1}=b_i>0$ for each $i=0,1,\ldots$ and $P_{i,i-j}=0$ for each $j>1$ and all $i$. In this case, Jiang \emph{et al$.$} \cite{jly14} have shown that 
\[
m_j(h)=\frac{1}{b_j\pi_j}\sum_{k=0}^j\widehat{h}(k)\pi_k=-\frac{1}{b_j\pi_j}\sum_{k=j+1}^\infty\widehat{h}(k)\pi_k;
\]     
see their equation (15), and note that the final equality follows from $\mathbb{E}\widehat{h}(\pi)=0$. This immediately gives the bound
\begin{equation}\label{eq:steinfactor}
\sup_{h\in\mathcal{H}}\sup_{j\in\mathbb{Z}^+}|m_j(h)|\leq\sup_{j\in\mathbb{Z}^+}\left\{\frac{\mathbb{P}(\pi>j)}{b_j\pi_j}\right\}.
\end{equation}

For a straightforward illustrative example, consider a simple random walk on $\mathbb{Z}^+$ with reflection at the origin, where $P_{0,0}=P_{i,i-1}=p$ for all $i\geq1$ and $P_{i,i+1}=1-p$ for all $i\geq0$, for some $p>1/2$ to ensure positive recurrence. We note that in this case $P$ is stochastically monotone, and that $\pi\sim\mbox{Geom}(\alpha)$ is geometrically distributed with parameter $\alpha=\frac{2p-1}{p}$ and mass function $\pi_k=\alpha(1-\alpha)^k$ for $k=0,1,\ldots$, so that \eqref{eq:steinfactor} gives $|m_j(h)|\leq1/(2p-1)$.

\subsection{Example: M/M/1 queue}\label{sec:mm1}

For our next example, we look beyond the class of birth-death processes into Markov chains of the type associated with GI/M/1 queues. For simplicity we will restrict our attention to the M/M/1 queue here. Although this may be formulated as a birth-death process, we instead use an alternative representation which makes it clear how this example can be extended to more general GI/M/1-type chains. 
Consider the M/M/1 queue in which customer interarrival times are exponentially distributed with mean $1/\lambda$ and service times are exponentially distributed with mean $1/\mu$, where we assume $\rho=\lambda/\mu<1$. Let $\{Z_t:t=0,1,\ldots\}$ be the Markov chain embedded at customer arrival times, which may be constructed as a stochastically monotone single-birth chain with stationary distribution $\pi\sim\mbox{Geom}(1-\rho)$; that is, $\pi_i=(1-\rho)\rho^i$ for $i=0,1,\ldots$. We define
\[
a_k=\frac{\rho}{1+\rho}\left(\frac{1}{1+\rho}\right)^k,
\] 
the probability that exactly $k$ customers are served in the time between two consecutive arrivals. We then define our transition matrix $P$ by writing $P_{i,j}=a_{i-j+1}$ for $j=1,2,\ldots,i+1$ and each $i$. The remaining non-zero entries of $P$ are those in the left-hand column, which are given by $P_{i,0}=\sum_{j>i}a_j=(1+\rho)^{-(i+1)}$, for $i\geq0$. Note that, as constructed, this is not a birth-death chain, and so we cannot use the bound \eqref{eq:steinfactor} here. We will need to calculate separately a bound on the increments $m_i(h)$ of the solution of \eqref{eq:poisson}. By \eqref{eq:jly}, these are given by $m_0(h)=(1+\rho)\rho^{-1}\widehat{h}(0)$ and
\[
m_i(h)=\frac{1+\rho}{\rho}\left(\widehat{h}(i)+\frac{1}{(1+\rho)^{i+1}}\sum_{k=0}^{i-1}(1+\rho)^km_k(h)\right),
\]
for $i=1,2,\ldots$. Solving this system of equations gives
\[
m_i(h)=\left(\frac{1+\rho}{\rho}\right)\widehat{h}(i)+\frac{1}{\rho^{i+1}}\sum_{j=0}^{i-1}\widehat{h}(j)\rho^j=\left(\frac{1+\rho}{\rho}\right)\widehat{h}(i)-\frac{1}{\rho^{i+1}}\sum_{j=i}^{\infty}\widehat{h}(j)\rho^j,
\]
where the final equality uses the fact that $\mathbb{E}\widehat{h}(\pi)=0$. For $h\in\mathcal{H}$ we therefore have
\begin{equation}\label{eq:steinfactor_mm1}
|m_i(h)|\leq\frac{1+\rho}{\rho}+\frac{1}{\rho^{i+1}}\sum_{j=i}^\infty\rho^j=\frac{2-\rho^2}{\rho(1-\rho)}.
\end{equation}

\section{Convergence to stationarity} \label{sec:convergence}

In the setting of Theorem \ref{thm:main}, we may choose $X$ to have the same distribution as $Z_t$ for some fixed $t\in\mathbb{Z}^+$. This lets us bound the total variation distance of $Z_t$ from stationarity. We set $Z_0=0$, and note that 
\[
\sum_{k=j+1}^\infty\sum_{i=0}^\infty\mathbb{P}(Z_t=i)P_{i,k}
=\sum_{k=j+1}^\infty\mathbb{P}(Z_{t+1}=k)
=\mathbb{P}(Z_{t+1}>j).
\]
Under the assumption that $P$ is stochastically monotone, we can couple our Markov chain in such a way that $\mathbb{P}(Z_{t+1}>j)\geq\mathbb{P}(Z_t>j)$ for all $j\in\mathbb{Z}^+$. Hence,
\[
\sum_{j=0}^\infty\left|\mathbb{P}(Z_t>j)-\mathbb{P}(Z_{t+1}>j)\right|=\mathbb{E}Z_{t+1}-\mathbb{E}Z_t,
\]
and Theorem \ref{thm:main} gives the following, in which we may choose a coupling of $Z_t$ and $Z_{t+1}$ to bound the expectation.
\begin{corollary}\label{cor:convergence}
Let $\{Z_t:t=0,1,\ldots\}$ be a positive recurrent and stochastically monotone single-birth Markov chain on $\mathbb{Z}^+$, as defined above, with transition matrix $P$ and stationary distribution $(\pi_0,\pi_1,\ldots)$. Then
\[
d_{TV}(Z_t,\pi)\leq\sup_{h\in\mathcal{H}}\sup_{l\in\mathbb{Z}^+}|m_l(h)|\left(\mathbb{E}Z_{t+1}-\mathbb{E}Z_t\right).
\]
\end{corollary}
We use the remainder of this section to illustrate the bound of Corollary \ref{cor:convergence} using the applications we introduced in Section \ref{sec:intro}.

\subsection{Example: Simple random walk with reflection}\label{sec:convergence_srw}

Consider the simple random walk with reflection at the origin that was introduced in Section \ref{sec:bd}, for which we know that $|m_j(h)|\leq1/(2p-1)$. To apply Corollary \ref{cor:convergence} it remains only to couple $Z_{t+1}$ and $Z_t$. To do this, we introduce a copy $\{Z^\prime_t:t=0,1,\ldots\}$ of $\{Z_t:t=0,1,\ldots\}$, with these two processes coupled as follows: with $Z_0=Z_0^\prime=0$, we let $Z_1^\prime$ be $0$ or $1$, with probability $p$ and $1-p$ respectively, independently of all else. The processes $\{Z_t:t=1,2,\ldots\}$ and $\{Z_{t+1}^\prime:t=1,2,\ldots\}$ then evolve using the same underlying sequence of independent Bernoulli trials so that, roughly speaking, one process moves in the positive direction at a given time if and only if the other process does also; the same is true of steps in the negative direction, except that we need to account for the reflection at the origin where a `step in the negative direction' corresponds to remaining at the origin. This continues until the first time $t$ at which $Z_t=Z_{t+1}^\prime$, following which the two processes move together. This happens at a time which is almost surely no greater than $T=\min\{t\geq1:Z_t^\prime=0\}$, the first return time to the origin, and we note that
\begin{multline*}
\mathbb{E}Z_{t+1}-\mathbb{E}Z_t=\mathbb{E}Z_{t+1}^\prime-\mathbb{E}Z_t\leq\mathbb{P}(T>t+1)\\
=\mathbb{P}(r^T\geq r^{t+2})\leq\frac{\mathbb{E}r^T}{r^{t+2}}=\frac{2pr+1-[1-4p(1-p)r^2]^{1/2}}{2r^{t+2}},
\end{multline*}  
for any $1\leq r\leq[4p(1-p)]^{-1/2}$; see Section XIV.4 of \cite{f60} for the final equality. Corollary \ref{cor:convergence} then gives us that, for any $1\leq r\leq[4p(1-p)]^{-1/2}$,
\[
d_{TV}(Z_t,\pi)\leq\left(\frac{2pr+1-[1-4p(1-p)r^2]^{1/2}}{2r^{2}(2p-1)}\right)r^{-t},
\]
yielding the expected geometric rate of convergence to stationarity (see also Example 7.1 of \cite{lt96}) and an explicit bound on the corresponding total variation distance.

\subsection{Example: M/M/1 queue}

Now consider the M/M/1 queue of Section \ref{sec:mm1}, for which the bound \eqref{eq:steinfactor_mm1} holds. To apply Corollary \ref{cor:convergence} we again only need to couple $Z_{t+1}$ and $Z_t$. We proceed similarly to above, introducing a coupled copy $\{Z_t^\prime:t=0,1,\ldots\}$ of $\{Z_t:t=0,1,\ldots\}$. With $Z_0=Z_0^\prime=0$, let $Z_1^\prime$ be 0 or 1, with probability $(1+\rho)^{-1}$ and $\rho(1+\rho)^{-1}$ respectively, independently of all else. We let subsequent arrivals occur at the same times in both processes, and hence $Z_t\leq Z_{t+1}^\prime$ almost surely for all $t$. As before, $\mathbb{E}Z_{t+1}-\mathbb{E}Z_t\leq\mathbb{P}(T>t+1)$, where $T$ is as in Section \ref{sec:convergence_srw}, since the numbers of customers in the two systems differ by at most one. We may write $T$ as $1+IN$, where $N$ is the number of customers served in a busy period of an M/M/1 queue initiated by the arrival of a single customer to an empty system, and $I$ is a Bernoulli random variable with mean $\rho(1+\rho)^{-1}$ independent of all else. We then have
\[
\mathbb{P}(T>t+1)=\frac{\rho}{1+\rho}\mathbb{P}(N>t)\leq\frac{\rho\mathbb{E}r^N}{(1+\rho)r^{t+1}}=\frac{1+\rho-[(1+\rho)^2-4\rho r]^{1/2}}{2(1+\rho)r^{t+1}},
\]  
for all $1\leq r\leq\frac{(1+\rho)^2}{4\rho}$; the expression for $\mathbb{E}r^N$ is well-known. Corollary \ref{cor:convergence} thus gives us, for any $1\leq r\leq\frac{(1+\rho)^2}{4\rho}$,
\[
d_{TV}(Z_t,\pi)\leq\left(\frac{(2-\rho^2)\left(1+\rho-[(1+\rho)^2-4\rho r]^{1/2}\right)}{2\rho(1-\rho)(1+\rho)r}\right)r^{-t}.
\] 

\section{Comparison of stationary distributions} \label{sec:comparison}

In this section we let $X$ have the stationary distribution of a positive recurrent Markov chain on $\mathbb{Z}^+$ with transition matrix $Q$ whose $(i,j)$th entry we denote by $Q_{i,j}$. We bound $d_{TV}(X,\pi)$ using Theorem \ref{thm:main}. In this setting we have
\[
\sum_{j=0}^\infty\left|\mathbb{P}(X>j)-\sum_{k=j+1}^\infty\sum_{i=0}^\infty\mathbb{P}(X=i)P_{i,k}\right|
\leq \sum_{j=0}^\infty\sum_{i=0}^\infty\mathbb{P}(X=i)\left|\sum_{k=j+1}^\infty\left(Q_{i,k}-P_{i,k}\right)\right|.
\]
As we illustrate in the examples below, results simplify further in the case where either $P$ dominates, or is dominated by, $Q$. That is, where either 
\begin{equation}\label{eq:dom}
\sum_{k>m}P_{i,k}\geq\sum_{k>m}Q_{i,k}
\end{equation}
for all $i,m\in\mathbb{Z}^+$, or the reverse inequality holds for all $i,m\in\mathbb{Z}^+$. In either of these cases we have
\[
\sum_{j=0}^\infty\sum_{i=0}^\infty\mathbb{P}(X=i)\left|\sum_{k=j+1}^\infty\left(Q_{i,k}-P_{i,k}\right)\right|
=\left|\mathbb{E}\sum_{j=0}^\infty\sum_{k=j+1}^\infty\left(Q_{X,k}-P_{X,k}\right)\right|
=\left|\mathbb{E}X-\mathbb{E}\sum_{j=0}^\infty\sum_{k=j+1}^\infty P_{X,k}\right|.
\]
We thus obtain the following from Theorem \ref{thm:main}.
\begin{corollary}\label{cor:comparison}
Let $P$ be the transition matrix of a positive recurrent single-birth Markov chain on $\mathbb{Z}^+$ with stationary distribution $(\pi_0,\pi_1,\ldots)$. Let $X$ have the stationary distribution of another positive recurrent Markov chain on $\mathbb{Z}^+$ with transition matrix $Q$. Then
\begin{equation}\label{eq:comparison}
d_{TV}(X,\pi)\leq\sup_{h\in\mathcal{H}}\sup_{l\in\mathbb{Z}^+}|m_l(h)|\sum_{i=0}^\infty\mathbb{P}(X=i)\sum_{j=0}^\infty\left|\sum_{k=j+1}^\infty\left(Q_{i,k}-P_{i,k}\right)\right|.
\end{equation} 
If, in addition, either \eqref{eq:dom} or the reverse inequality holds for all $i,m\in\mathbb{Z}^+$ then
\begin{equation}\label{eq:comparison2}
d_{TV}(X,\pi)\leq\sup_{h\in\mathcal{H}}\sup_{l\in\mathbb{Z}^+}|m_l(h)|\left|\mathbb{E}X-\mathbb{E}\sum_{j=0}^\infty\sum_{k=j+1}^\infty P_{X,k}\right|.
\end{equation}
\end{corollary}

\subsection{Example: Birth--death chains}

In the setting where both $P$ and $Q$ are transition matrices of birth--death chains, combining \eqref{eq:steinfactor} and \eqref{eq:comparison} gives the following bound between the corresponding stationary distributions:
\[
d_{TV}(X,\pi)\leq\sup_{j\in\mathbb{Z}^+}\left\{\frac{\mathbb{P}(\pi>j)}{P_{j,j+1}\pi_j}\right\}\bigg(\mathbb{E}|P_{X,X-1}-Q_{X,X-1}|+\mathbb{E}|P_{X,X+1}-Q_{X,X+1}|\bigg).
\]
A similar bound applies in other settings, for example in approximating the stationary distribution of a single-birth chain by that of a birth--death chain.

\subsection{Example: Geometric approximation}

We give two applications of \eqref{eq:comparison2} to approximation by a geometric distribution for the stationary distribution associated with our transition matrix $Q$, using transition matrices $P$ which have a geometric stationary distribution.

\subsubsection{Simple random walk with reflection}

Let $P$ be the transition matrix of the simple random walk with reflection, as in Section \ref{sec:bd}. For $i\geq0$,
\[
\sum_{k>m}P_{i,k}=\left\{\begin{array}{ll} 1 & \text{if }m<i-1,\\ 1-p, & \text{if }i-1\leq m\leq i,\\0, & \text{if }m>i.\end{array}\right.
\]
If we assume that either \eqref{eq:dom} or the reverse inequality holds for all $i,m\in\mathbb{Z}^+$, we may apply \eqref{eq:comparison2} to give an easily computed bound in the approximation of the stationary distribution associated with $Q$ by $\pi\sim\mbox{Geom}(\alpha)$, where $\alpha=(2p-1)/p$. Recalling that $|m_j(h)|\leq1/(2p-1)$ in this setting, a simple calculation shows that \eqref{eq:comparison2} gives
\[
d_{TV}(\pi,X)\leq\left|1-\frac{p}{2p-1}\mathbb{P}(X=0)\right|,
\]
which we note is typically easy to evaluate and gives zero in the case where $X\sim\mbox{Geom}(\alpha)$.

\subsubsection{M/M/1 queue}

With $P$ as in Section \ref{sec:mm1} we may use the bound \eqref{eq:steinfactor_mm1} on $m_l(h)$. We further note that $\sum_{k>m}P_{i,k}=1-(1+\rho)^{m-i-1}$
for $i\geq m\geq0$, and $\sum_{k>m}P_{i,k}=0$ for $i<m$. 
We may therefore apply \eqref{eq:comparison2} for a Markov chain with transition matrix $Q$ satisfying $\sum_{k=0}^m Q_{i,k}\leq(1+\rho)^{m-i-1}$ for all $i\geq m\geq 0$. That is, if $Q$ is such that the total mass in the first $m$ elements of row $i$ decreases geometrically in $i$, and increases at most geometrically in $m$, Corollary \ref{cor:comparison} gives us a bound in the approximation of the corresponding stationary distribution by a geometric distribution with an appropriately chosen parameter $1-\rho$. Specifically, \eqref{eq:comparison2} yields
\[
d_{TV}(\pi,X)\leq\frac{2-\rho^2}{\rho(1-\rho)}\left\{\frac{1-\rho}{\rho}-\frac{1}{\rho}\mathbb{E}\left(\frac{1}{1+\rho}\right)^{X+1}\right\}.
\]
Note that this upper bound is zero if $X\sim\mbox{Geom}(1-\rho)$, as expected.

\subsection{Example: Truncation} \label{sec:truncation}

We conclude with a final application, to the truncation of the state space of our single-birth Markov chain with transition matrix $P$. Let $Q=Q_{(n)}$ denote the $(n+1)\times(n+1)$ northwest truncation of $P$, augmented to be a valid transition matrix by replacing $P_{n,j}$ in the final row of $P$ by $Q_{n,j}=P_{n,j}+\nu_jP_{n,n+1}$ for each $j\geq0$, where $(\nu_0,\nu_1,\ldots,\nu_n)$ is a probability distribution, and we denote by $\nu$ a random variable with this distribution. Let $X$ have the stationary distribution associated with the transition matrix $Q$; this is a special case of the truncation problem for discrete-time Markov chains studied by many authors. Tweedie \cite{t98} shows that the corresponding stationary probabilities converge to $\pi_j$ in the case of a geometrically ergodic chain, a stochastically monotone chain or one dominated by a stochastically monotone chain, when the augmentation is in the first or last column only. These results have since been generalised (see, for example, \cite{igl22,l10,ll18} and references therein) and given improved error bounds. Our purpose here is not to compete with these general bounds, but to illustrate the straightforward application of Corollary \ref{cor:comparison} in this setting and to note the explicit bound it yields.

Since $\nu_j\geq0$ for each $j\leq n$, it is clear that $\sum_{k>m}Q_{i,k}\leq\sum_{k>m}P_{i,k}$ for each $i$ and $m$. Since
\begin{align*}
\mathbb{E}\sum_{j=0}^\infty\sum_{k=j+1}^\infty\left(P_{X,k}-Q_{X,k}\right)&=\mathbb{P}(X=n)\sum_{j=0}^n\left(P_{n,n+1}-\sum_{k=j+1}^n\nu_kP_{n,n+1}\right)\\
&=\mathbb{P}(X=n)\left(n+1-\mathbb{E}\nu\right)P_{n,n+1},
\end{align*}
Corollary \ref{cor:comparison} gives
\begin{equation}\label{eq:truncation}
d_{TV}(X,\pi)\leq\sup_{h\in\mathcal{H}}\sup_{l\in\mathbb{Z}^+}|m_l(h)|\mathbb{P}(X=n)\left(n+1-\mathbb{E}\nu\right)P_{n,n+1}.
\end{equation}

As a simple illustrative example, suppose that $P$ is the transition matrix of a stochastically monotone birth--death chain with stationary distribution given by $\pi_k=\alpha(1-\alpha)^k$ for some $\alpha\in(0,1)$ and all $k=0,1,\ldots$. The simple random walk with reflection of Section \ref{sec:bd} is an example of such a chain. Stochastic monotonicity of $P$ gives us that $\mathbb{P}(X=n)\leq\mathbb{P}(\pi\geq n)=(1-\alpha)^n$, and \eqref{eq:truncation} becomes  
\[
d_{TV}(X,\pi)\leq\frac{(n+1-\mathbb{E}\nu)P_{n,n+1}}{\alpha\inf_{j\in\mathbb{Z}^+}\{P_{j,j+1}\}}(1-\alpha)^{n+1},
\]
giving the same rate as the lower bound $d_{TV}(X,\pi)\geq\mathbb{P}(\pi>n)=(1-\alpha)^{n+1}$ if $P_{j,j+1}$ is bounded away from zero.

\end{document}